\documentclass[12pt]{article}
\usepackage{latexsym,amssymb,amsmath}
\begin{document}

\baselineskip=20pt

\newcommand{\la}{\langle}
\newcommand{\ra}{\rangle}
\newcommand{\psp}{\vspace{0.4cm}}
\newcommand{\pse}{\vspace{0.2cm}}
\newcommand{\ptl}{\partial}
\newcommand{\dlt}{\delta}
\newcommand{\sgm}{\sigma}
\newcommand{\al}{\alpha}
\newcommand{\be}{\beta}
\newcommand{\G}{\Gamma}
\newcommand{\gm}{\gamma}
\newcommand{\vs}{\varsigma}
\newcommand{\lmd}{\lambda}
\newcommand{\td}{\tilde}
\newcommand{\vf}{\varphi}
\newcommand{\rd}{\mbox{Rad}}
\newcommand{\ad}{\mbox{ad}}
\newcommand{\stl}{\stackrel}
\newcommand{\ol}{\overline}
\newcommand{\ul}{\underline}
\newcommand{\es}{\epsilon}
\newcommand{\dmd}{\diamond}
\newcommand{\clt}{\clubsuit}
\newcommand{\vt}{\vartheta}
\newcommand{\ves}{\varepsilon}
\newcommand{\kn}{\mbox{ker}}
\newcommand{\for}{\mbox{for}}
\newcommand{\rad}{\mbox{Rad}}
\newcommand{\ups}{\Upsilon}

\begin{center}{\Large \bf Derivation-Simple Algebras and the Structures}\end{center}
\begin{center}{\Large \bf  of Generalized  Lie Algebras of Witt Type}
\end{center}
\vspace{0.2cm}

\begin{center}{\large Yucai Su$^{\ast}$, Xiaoping Xu$^{\dag}$ and Hechun Zhang$^{\ddag}$ }\end{center}

* Department of Applied Mathematics, Shanghai Jiaotong University, 1954 Huashan Road, Shanghai 200030, P. R. China.

\dag Department of Mathematics, The Hong Kong University of Science \& Technology, Clear Water Bay, Kowloon, Hong Kong, P. R. China.

\ddag  Department of Applied Mathematics, Tsinghua University, Beijing 100084, P. R. China.

\vspace{0.3cm}

\begin{center}{\Large \bf Abstract}\end{center}
\vspace{0.2cm}

{\small We classify all the pairs of a commutative associative algebra with an identity element and its finite-dimensional commutative locally-finite derivation subalgebra such that the commutative associative algebra is derivation-simple with respect to the derivation subalgebra over an algebraically closed field with characteristic 0. Such pairs are the fundamental ingredients for constructing generalized simple Lie algebras of Cartan type. Moreover, we determine the isomorphic classes of the generalized simple Lie algebras of Witt Type. The structure space of these algebras is given explicitly.}

\section{Introduction}

Simple Lie algebras of Cartan type are important geometrically-natural infinite-dimensional Lie algebras in mathematics. The fundamental ingredients of constructing Lie algebras of Cartan type are the pairs of a polynomial algebra and the derivation subalgebra of the operators of taking partial derivatives. The abstract definition of generalized Lie algebras of  Cartan type by derivations  appeared in Kac's work [Ka1]. However, it is still a question of how to construct new explicit generalized simple Lie algebras of Cartan type. Kawamoto [K] constructed new generalized simple Lie algebras of Witt type by the pairs of the group algebra of an additive subgroup of $\Bbb{F}^n$ and the derivation subalgebra of the grading operators, where $n$ is a positive integer and $\Bbb{F}$ is a field with characteristic 0. One can view the operators of taking partial derivatives of  the polynomial algebra in several variables as {\it down-grading operators}. Using the pairs of the tensor algebra of the group algebra of the direct sum of finite number of additive subgroups of $\Bbb{F}$ with the polynomial algebra in several variables, and the derivation subalgebra of the grading operators and down-grading operators, Osborn [O] constructed new generalized simple Lie algebras of Cartan type. In [DZ1], the authors generalized Kawamoto's work by picking out certain subalgebras. Their construction is also equivalent to generalizing  Osborn's Lie algebras of Witt type by adding certain  diagonal elements of $\Bbb{F}^n$ into the group. 

Passman [P] proved that the generalized simple Lie algebras of Witt type constructed from the pairs of a commutative associative algebra with an identity element and its commutative derivation subalgebra are simple Lie algebras if and only if the commutative associative algebra is derivation-simple with respect to the derivation subalgebra. In [X1], the second author of this paper constructed new explicit generalized simple Lie algebras of Cartan type, based on the pairs of the tensor algebra of the group algebra of an additive subgroup of $\Bbb{F}^n$ with the polynomial algebra in several variables and the derivation subalgebra of the mixtures of the grading and down-grading operators. The algebras in [X1] are the most general known explicit examples of generalized simple Lie algebras of Cartan type.  A natural question is how far it is from the generalized simple Lie algebras of Witt type in [X1] to those abstractly determined by Passman [P]. In this paper, we shall show that the generalized  simple Lie algebras of Witt type determined in [P] are essentially the same as those explicitly constructed in [X1] under certain locally-finite conditions. It seems to us that further constructing new generalized simple Lie algebras of Cartan type that are essentially different from those in [X1] and some extensions in Chapter 6 of [X2] is extremely difficult. Of course, one can get some new generalized simple Lie algebras of Cartan type through replacing the commutative associative algebra used in [X1] by a certain subalgebra of its topological completion. However, such a construction is not essential from algebraic point of view.

     Zhao determined the isomorphic classes of the generalized simple Lie algebras of Witt type found in [DZ1]. In this paper, we shall determine the isomorphic classes of the generalized simple Lie algebras of Witt type constructed in [X1], which are more general than those in [DZ1]. The structure space of the generalized simple Lie algebras of Witt type constructed in [X1] will be given explicitly.

      Below, we shall give a more detailed description of our results. Throughout this paper, $\Bbb{F}$ denotes an algebraically closed field with characteristic 0. All the vector spaces (algebras) without specifying field  are assumed over $\Bbb{F}$. We always assume that an associative algebra has an identity element. Moreover, we denote by $\Bbb{Z}$ the ring of integers and by $\Bbb{N}$ the set of nonnegative integers.

      Let ${\cal A}$ be a commutative associative algebra. A {\it derivation} $d$ of ${\cal A}$ is a linear transformation on ${\cal A}$ such that
$$d(uv)=d(u)v+ud(v)\qquad\for\;u,v\in{\cal A}.\eqno(1.1)$$
The space $\mbox{Der}\:{\cal A}$ of derivations forms a Lie algebra with respect to the operation:
$$[d_1,d_2]=d_1 d_2-d_2d_1\qquad\for\;\;d_1,d_2\in \mbox{Der}\:{\cal A}.\eqno(1.2)$$
For  $u\in{\cal A}$ and $d\in \mbox{Der}\:{\cal A}$, we define
$$(ud)(v)=ud(v)\qquad\for\;\;v\in{\cal A}.\eqno(1.3)$$
Since ${\cal A}$ is a commutative associative algebra, $ud$ is also a derivation. In particular, $\mbox{Der}\:{\cal A}$ is an ${\cal A}$-module.

A linear transformation $T$ on a vector space $V$ is called {\it locally-finite} if 
$$\mbox{dim span}\:\{T^m(v)\mid m\in\Bbb{N}\}<\infty\eqno(1.4)$$
for any $v\in V$. A set of linear transformations is called {\it locally-finite} if all its elements are locally-finite. 

For a commutative associative algebra ${\cal A}$ and a derivation subspace ${\cal D}$, ${\cal A}$ is called {\it derivation-simple} with respect to ${\cal D}$ if there does not exist a subspace ${\cal I}$ of ${\cal A}$ such that ${\cal I}\neq\{0\},{\cal A}$ and
$$u{\cal I},\;d({\cal I})\subset {\cal I}\qquad\for\;\;u\in{\cal A},\;d\in{\cal D}.\eqno(1.5)$$
Moreover, the derivation subspace ${\cal D}$ is a {\it commutative subalgebra} if
$$d_1d_2=d_2d_1\qquad\for\;\;d_1,d_2\in \mbox{Der}\:{\cal A}.\eqno(1.6)$$

In this paper, we shall first give in Section 2 a complete classification of the pairs of a  commutative associative algebra ${\cal A}$ and a finite-dimensional locally-finite commutative derivation subalgebra ${\cal D}$ such that ${\cal A}$ is derivation-simple with respect to ${\cal D}$. For such a pair $({\cal A},{\cal D})$ with ${\cal D}\neq\{0\}$, Passman's Theorem [P] tells us that
$${\cal W}={\cal A}{\cal D}\eqno(1.7)$$
forms a simple Lie algebra, which is called a {\it generalized Lie algebra of Witt type}. In Section 3, we shall determine the isomorphic classes of the generalized simple Lie algebras of the form (1.7). The structure space will be presented explicitly.  From purely Lie algebra structure point of view, it is enough to consider the Lie algebras of the form (1.7) with
$$\bigcap_{d\in{\cal D}}\kn_{d}=\Bbb{F}.\eqno(1.8)$$

\section{Derivation-Simple Algebra}

In this section, we shall classify the pairs of a commutative associative algebra ${\cal A}$ and a finite-dimensional locally-finite commutative derivation subalgebra ${\cal D}$ such that ${\cal A}$ is derivation-simple with respect to ${\cal D}$.

We start with constructions of such pairs. Let $k_1$ and $k_2$ be two nonnegative integers such that
$$k=k_1+k_2>0.\eqno(2.1)$$
Let $\G$ be an additive subgroup of $\Bbb{F}^k$ and let $\Bbb{F}_1$ be an extension field of $\Bbb{F}$. Suppose that $f(\cdot,\cdot):\G\times \G\rightarrow \Bbb{F}^{\times}=\Bbb{F}\setminus\{0\}$ is a map such that 
$$f(\al,\be)f(\al+\be,\gm)=f(\al,\be+\gm)f(\be,\gm),\;\;f(\al,\be)=f(\be,\al),\;\;f(\al,0)=1\eqno(2.2)$$
for $\al,\be,\gm\in \G$. 
Denote by $\Bbb{F}_1[t_1,t_2,...,t_{k_1}]$  the algebra of polynomials in $k_1$ variables over $\Bbb{F}_1$. Let ${\cal A}(k_1,k_2;\G,\Bbb{F}_1,f)$ be a free $\Bbb{F}_1[t_1,t_2,...,t_{k_1}]$-module with the basis
$$\{x^{\al}\mid \al\in\G\}.\eqno(2.3)$$
Viewing ${\cal A}(k_1,k_2;\G,\Bbb{F}_1,f)$ as a vector space over $\Bbb{F}$, we define an algebraic operation ``$\cdot$'' on ${\cal A}(k_1,k_2;\G,\Bbb{F}_1,f)$ by
$$(\zeta x^{\al})\cdot (\eta x^{\be})=f(\al,\be)\zeta\eta x^{\al+\be}\qquad\for\;\;\zeta,\eta\in \Bbb{F}_1[t_1,t_2,...,t_{k_1}],\;\al,\be\in\G.\eqno(2.4)$$
Then $({\cal A}(k_1,k_2;\G,\Bbb{F}_1,f),\cdot)$ forms a commutative associative algebra over $\Bbb{F}$ with $x^{0}$ as the identity element, which is denoted as $1$ for convenience. We refer Subsection 5.4.2 of [X2] for the details of this algebra. When the context is clear, we shall omit the notion ``$\cdot$'' in any associative algebra product.

We define the linear transformations 
$$\{\ptl_{t_1},...,\ptl_{t_{k_1}},\ptl^{\ast}_1,...,\ptl^{\ast}_k\}\eqno(2.5)$$
on ${\cal A}(k_1,k_2;\G,\Bbb{F}_1,f)$ by 
$$\ptl_{t_i}(\zeta x^{\al})=\ptl_{t_i}(\zeta)x^{\al},\;\;\ptl^{\ast}_j(\zeta x^{\al})=\al_j \zeta x^{\al}\eqno(2.6)$$
for $\zeta\in \Bbb{F}_1[t_1,t_2,...,t_{k_1}],\;\al=(\al_1,...,\al_k)\in\G,$ where on $\Bbb{F}_1[t_1,t_2,...,t_{k_1}]$,  $\ptl_{t_i}$ are the operators of taking partial derivative with respect to $t_i$ over $\Bbb{F}_1$. Then $\{\ptl_{t_1},...,\ptl_{t_{k_1}},\ptl^{\ast}_1,...,\ptl^{\ast}_k\}$ are mutually commutative derivations of ${\cal A}(k_1,k_2;\G,\Bbb{F}_1,f)$. The derivations  $\{\ptl_{t_1},...,\ptl_{t_{k_1}}\}$ are called {\it down-grading operators} and $\{\ptl^{\ast}_1,...,\ptl^{\ast}_k\}$ are called {\it grading operators} of the algebra ${\cal A}(k_1,k_2;\G,\Bbb{F}_1,f)$. Given $m,n\in\Bbb{Z}$ with $m<n$, we shall use the following notion 
$$\ol{m,n}=\{m,m+1,m+2,...,n\}\eqno(2.7)$$
throughout this paper. We also treat $\ol{m,n}=\emptyset$ when $m>n$. 

Choose
$$\{\bar{\ptl}_{k_1+1},...,\bar{\ptl}_k\}\subset \sum_{j=1}^{k_1}\Bbb{F}_1\ptl_{t_j}\eqno(2.8)$$
(cf. (1.3)). We set
$$\ptl_i=\ptl_i^{\ast}+\ptl_{t_i},\;\;\ptl_{k_1+j}=\ptl^{\ast}_{k_1+j}+\bar{\ptl}_j\qquad\for\;\;i\in\ol{1,k_1},\;j\in\ol{1,k_2}.\eqno(2.9)$$
Then $\{\ptl_i\mid i\in\ol{1,k}\}$ is a  set of derivations. Set
$${\cal D}=\sum_{i=1}^k\Bbb{F}\ptl_i.\eqno(2.10)$$
Note that ${\cal D}$ is a finite-dimensional locally-finite commutative derivation subalgebra.
\psp

{\bf Theorem 2.1}. {\it Let} ${\cal A}$ {\it be a commutative associative algebra  and let} ${\cal D}$ {\it be a finite-dimensional locally-finite commutative derivation subalgebra}. {\it The algebra} ${\cal A}$ {\it is derivation-simple with respect to} ${\cal D}$ {\it if and only if the algebra} ${\cal A}$ {\it is isomorphic to the algebra of the form} ${\cal A}(k_1,k_2;\G,\Bbb{F}_1,f)$ {\it and the derivation subalgebra} ${\cal D}$ {\it  of the form (2.10)}.
\psp

{\it Proof}. Let us first prove that the algebra ${\cal A}(k_1,k_2;\G,\Bbb{F}_1,f)$ is derivation-simple with respect to the derivation subalgebra ${\cal D}$ in (2.10). Let ${\cal I}$ be a nonzero ${\cal D}$-invariant ideal of ${\cal A}(k_1,k_2;\G,\Bbb{F}_1,f)$. For any $\al=(\al_1,...,\al_k)\in\G$, we define
$$\bar{\cal A}_{\al}=\Bbb{F}_1[t_1,t_2,...,t_{k_1}]x^{\al}.\eqno(2.11)$$
Then 
$$\bar{\cal A}_{\al}=\{u\in {\cal A}(k_1,k_2;\G,\Bbb{F}_1,f)\mid (\ptl_i-\al_i)^m(u)=0\;\mbox{for some}\;m\in\Bbb{N};\;i\in\ol{1,k}\}\eqno(2.12)$$
and 
$${\cal A}(k_1,k_2;\G,\Bbb{F}_1,f)=\bigoplus_{\al\in\G}\bar{\cal A}_{\al}.\eqno(2.13)$$
Thus we have
$${\cal I}=\bigoplus_{\al\in\G}{\cal I}_{\al},\qquad{\cal I}_{\al}=\bar{\cal A}_{\al}\bigcap {\cal I}.\eqno(2.14)$$
If ${\cal I}_{\be}\neq \{0\}$ for some $\be\in\G$, then $\{0\}\neq x^{-\be}({\cal I}_{\be})\subset {\cal I}_{0}$. Thus we always have ${\cal I}_{0}\neq\{0\}$. Since
$$\ptl_i|_{{\cal I}_{0}}=\ptl_{t_i}|_{{\cal I}_{0}}\qquad\for\;\;j\in\ol{1,k_1}\eqno(2.15)$$
by (2.9), we have $1\in {\cal I}$. Hence ${\cal A}(k_1,k_2;\G,\Bbb{F}_1,f)\cdot 1\subset {\cal I}$. So ${\cal I}={\cal A}(k_1,k_2;\G,\Bbb{F}_1,f)$. That is, ${\cal A}(k_1,k_2;\G,\Bbb{F}_1,f)$ is derivation-simple with respect to ${\cal D}$ in (2.10).

Next we assume that ${\cal A}$ is a commutative associative algebra and  ${\cal D}$ is a finite-dimensional locally-finite commutative derivation subalgebra such that ${\cal A}$  is derivation-simple with respect to ${\cal D}$. Denote by ${\cal D}^{\ast}$ the linear functions from ${\cal D}$ to $\Bbb{F}$, which forms a vector space with respect to the addition and scalar multiplication of functions.
 Since $\Bbb{F}$ is algebraically closed and ${\cal D}$ is finite-dimensional, commutative and locally-finite, we have
$${\cal A}=\bigoplus_{\al\in{\cal D}^{\ast}}{\cal A}_{\al},\;{\cal A}_{\al}=\{u\in{\cal A}\mid (d-\al(d))^m(u)=0\;\;\mbox{for}\;\;d\in{\cal D}\;\mbox{and some}\;m\in\Bbb{N}\}.\eqno(2.16)$$
Denote
$$\G=\{\al\mid {\cal A}_{\al}\neq\{0\}\}.\eqno(2.17)$$
For any $\al\in{\cal D}^{\ast}$ and $n\in\Bbb{N}$, we define
$${\cal A}_{\al}^{(n)}=\{u\in{\cal A}\mid (d-\al(d))^{n+1}(u)=0\;\mbox{for}\;d\in{\cal D}\}.\eqno(2.18)$$
Then 
$${\cal A}_{\al}=\bigcup_{n=0}^{\infty}{\cal A}_{\al}^{(n)}\eqno(2.19)$$
and
$${\cal A}_{\al}=\{0\}\Longleftrightarrow {\cal A}_{\al}^{(0)}=\{0\}.\eqno(2.20)$$
We call a nonzero element in ${\cal A}_{\al}^{(0)}$ a {\it root vector}.

For any root vector $u$, ${\cal A}u$ is a ${\cal D}$-invariant ideal of ${\cal A}$. Thus ${\cal A}u={\cal A}$. In particular, $vu=1_{\cal A}$ for some $ v\in{\cal A}$. So any root vector is always invertible.
For $u\in{\cal A}_{\al}^{(0)}$ with $\al\in\G$ and $d\in{\cal D}$, we have
$$0=d(1)=d(uu^{-1})=d(u)u^{-1}+ud(u^{-1})=\al(d)uu^{-1}+ud(u^{-1})=\al(d)+ud(u^{-1}),\eqno(2.21)$$
which implies 
$$d(u^{-1})=ud(u^{-1})u^{-1}=-\al(d)u^{-1}.\eqno(2.22)$$
Hence $-\al\in\G$.
By (1.1), we have
$${\cal A}_{\al}^{(0)}\cdot {\cal A}_{\be}^{(0)}\subset {\cal A}_{\al+\be}^{(0)}\qquad\for\;\;\al,\be\in\G.\eqno(2.23)$$
Expression (2.23) and the invertibility of root vectors implies
$${\cal A}_{\al}^{(0)}\cdot {\cal A}_{\be}^{(0)}={\cal A}_{\al+\be}^{(0)}\qquad\for\;\;\al,\be\in\G.\eqno(2.24)$$
In particular, we obtain
$$\al+\be\in\G\qquad\for\;\;\al,\be\in\G.\eqno(2.25)$$
Thus $\G$ is an additive subgroup of ${\cal D}^{\ast}$.

Set
$$\Bbb{F}_1={\cal A}^{(0)}_0.\eqno(2.26)$$
Then $\Bbb{F}_1$ is an extension field of $\Bbb{F}$ by the invertibility of root vectors and (2.24). Suppose ${\cal A}_0\neq \Bbb{F}_1$. Note that for $u\in{\cal A}^{(m)}_{\al}$ and $v\in{\cal A}^{(n)}_{\be}$, we have
$$(d-(\al+\be)(d))^{m+n+1}(uv)=\sum_{i=0}^{m+n+1}(^{m+n+1}_{\;\;\;\;\;\:i})(d-\al(d))^i(u)(d-\be(d))^{m+n+1-i}(v)=0,\eqno(2.27)$$
that is, $uv\in {\cal A}_{\al+\be}^{(m+n)}$. So
$${\cal A}^{(m)}_{\al}\cdot {\cal A}^{(n)}_{\be}\subset {\cal A}_{\al+\be}^{(m+n)}\qquad\for\;\;\al,\be\in\G,\;m,n\in\Bbb{N}.\eqno(2.28)$$
In particular, 
$$\Bbb{F}_1{\cal A}^{(m)}_{\al}={\cal A}^{(m)}_{\al}\qquad\for\;\;\al\in\G,\;m\in\Bbb{N}.\eqno(2.29)$$
Hence each ${\cal A}^{(m)}_{\al}$ is a vector space over $\Bbb{F}_1$. 

For any $v\in{\cal A}_0^{(1)}$, ${\cal D}(v)\subset \Bbb{F}_1$ and
$${\cal D}(v)=\{0\}\Longleftrightarrow v\in \Bbb{F}_1.\eqno(2.30)$$ 
Set
$$H=\Bbb{F}_1{\cal D}\eqno(2.31)$$
(cf. (1.3)). Expression (2.30) implies that ${\cal A}^{(1)}_0/\Bbb{F}_1$ is isomorphic to a subspace of the space $\mbox{Hom}_{\Bbb{F}_1}(H,\Bbb{F}_1)$ over $\Bbb{F}_1$. 
By linear algebra, there exist subsets 
$$\{\ptl_1,...,\ptl_{k_1}\}\subset {\cal D},\;\;\{t_1,t_2,...,t_{k_1}\}\subset {\cal A}^{(1)}_0\eqno(2.32)$$
such that 
$${\cal A}^{(1)}_0=\Bbb{F}_1+\sum_{i=1}^{k_1}\Bbb{F}_1t_i,\;\;\ptl_i(t_j)=\dlt_{i,j}\qquad\for\;\;i,j\in\ol{1,k_1}.\eqno(2.33)$$
We write
$$H_1=\sum_{i=1}^{k_1}\Bbb{F}_1\ptl_i,\;\;H_2=\{d\in H\mid d({\cal A}^{(1)}_0)=\{0\}\}.\eqno(2.34)$$
Then we have
$$H=H_1\oplus H_2.\eqno(2.35)$$

Set
$$\bar{\cal A}_0=\sum_{n_i\in\Bbb{N};i\in\ol{1,k_1}}\Bbb{F}_1t_1^{n_1}t_2^{n_2}\cdots t_{k_1}^{n_{k_1}}\subset{\cal A}_0.\eqno(2.36)$$
Then $\bar{\cal A}_0$ forms a subalgebra of ${\cal A}$ and is isomorphic to $\Bbb{F}_1[t_1,...,t_{k_1}]$ when we view $t_i$ as variables by the second equation in (2.33). Moreover, ${\cal A}^{(1)}_0\subset \bar{\cal A}_0$ by the first equation in (2.33). Suppose ${\cal A}^{(m)}_0\subset \bar{\cal A}_0$ for some $1\leq m\in\Bbb{N}$. Note that
$$H_2(\bar{\cal A}_0)=\{0\}\eqno(2.37)$$
by (1.1) and (2.34).

By (2.18),
$$d({\cal A}_0^{(m+1)})\subset {\cal A}^{(m)}_0\subset\bar{\cal A}_0\qquad\for\;\;d\in H.\eqno(2.38)$$
For any $u\in {\cal A}_0^{(m+1)}$, there exists $u_1\in\bar{\cal A}_0$ such that
$$\ptl_1(u)=\ptl_1(u_1)\eqno(2.39)$$
by the derivation property of a polynomial algebra.
Similarly, we can find $u_2,...,u_{k_1}\in\bar{\cal A}_0$ such that
$$\ptl_i(u-\sum_{j=1}^iu_j)=0,\;\;\ptl_1(u_i)=\cdots =\ptl_{i-1}(u_i)=0\qquad\for\;\;i\in\ol{2,k_1}\eqno(2.40)$$
by induction on $i$. Thus we have
$$\ptl_i(u-\sum_{j=1}^{k_1}u_j)=0\qquad\mbox{for}\;\;i\in\ol{1,k_1}.\eqno(2.41)$$
For any $d\in H_2$, we have
$$d^2(u-\sum_{j=1}^{k_1}u_j)\in d({\cal A}^{(m)}_0)\subset d(\bar{\cal A}_0)=\{0\}\eqno(2.42)$$
by (2.37) and (2.38). Hence
$$u-\sum_{j=1}^{k_1}u_j\in{\cal A}_0^{(1)}\eqno(2.43)$$
by (2.18), (2.41) and (2.42). By the second equation (2.34), we get
$$d(u-\sum_{j=1}^{k_1}u_j)=0\qquad\for\;\;d\in H_2.\eqno(2.44)$$
Expressions (2.41) and (2.44) imply
$$u-\sum_{j=1}^{k_1}u_j\in\Bbb{F}_1,\eqno(2.45)$$
that is, $u\in\bar{\cal A}_0$. Therefore, ${\cal A}_0^{(m+1)}\subset \bar{\cal A}_0$. By induction on $m$, we obtain
$${\cal A}_0=\bar{\cal A}_0.\eqno(2.46)$$

The case  ${\cal A}_0=\Bbb{F}_1$ can be viewed as in the general case ${\cal A}_0=\Bbb{F}_1[t_1,...,t_{k_1}]$ with $k_1=0$.

For any $\al\in\G$, we take  $0\neq u\in{\cal A}^{(0)}_{\al}$ and have
$$u^{-1}{\cal A}_{\al}\subset {\cal A}_0,\;\;u{\cal A}_0\subset {\cal A}_{\al}\eqno(2.47)$$
by (2.22) and (2.28). Hence
$$u{\cal A}_0={\cal A}_{\al}.\eqno(2.48)$$
In particular, we have
$${\cal A}_{\al}^{(0)}=\Bbb{F}_1u\eqno(2.49)$$
is one-dimensional over $\Bbb{F}_1$. Choose 
$$x^0=1,\;\;0\neq x^{\al}\in{\cal A}_{\al}^{(0)}\qquad\mbox{for}\;\;0\neq \al\in\G.\eqno(2.50)$$
By (2.24) and (2.49), we have
$$x^{\al}x^{\be}=f(\al,\be)x^{\al+\be}\qquad\mbox{with}\;\;f(\al,\be)\in\Bbb{F}_1\;\;\mbox{for}\;\;\al,\be\in\G.\eqno(2.51)$$
Moreover, the first equation (2.50), the commutativity and associativity of ${\cal A}$ imply (2.2). Furthermore, (2.48) and (2.49) imply
$${\cal A}=\sum_{\al\in\G}{\cal A}_0x^{\al}.\eqno(2.52)$$

Observe that $\{\ptl_1,...,\ptl_{k_1}\}$ is an $\Bbb{F}$-linearly independent subset of ${\cal D}$ by (2.33). Extend it to an $\Bbb{F}$-basis $\{\ptl_1,...,\ptl_k\}$ of ${\cal D}$. Identifying
$$\al\longleftrightarrow (\al(\ptl_1),...,\al(\ptl_k))\qquad\for\;\;\al\in\G,\eqno(2.53)$$
we can view $\G$ as an additive subgroup of $\Bbb{F}^k$. Moreover, by  (2.31)-(2.37) and (2.46), we have
$$\ptl_j|_{{\cal A}_0}\in (\sum_{i=1}^{k_1}\Bbb{F}_1\ptl_i)|_{{\cal A}_0}\qquad\mbox{for}\;\;j\in\ol{k_1+1,k}.\eqno(2.54)$$
Therefore, the algebra ${\cal A}$ is isomorphic to the algebra ${\cal A}(k_1,k_2;\G,\Bbb{F}_1,f)$ with $k_2=k-k_1$ and ${\cal D}$ is of the form (2.10) by (2.52). This completes the proof of Theorem 2.1.

\section{Generalized Simple Lie Algebras of Witt Type}

In this section, we shall determine the structure space of the generalized simple Lie algebras of Witt type constructed in [X1]; namely, the isomorphic classes of the Lie algebra ${\cal W}$ of the form (3.1) with 
$$\Bbb{F}_1=\Bbb{F}\eqno(3.2)$$
for different $k_1, k_2$ and $\G$. First we need to rewrite ${\cal W}$ in a more compact form. Since $\Bbb{F}$ is algebraically closed, the algebra ${\cal A}(k_1,k_2;\G,\Bbb{F},f)$ is isomorphic to the semi-group algebra ${\cal A}(k_1,k_2;\G,\Bbb{F},{\bf 1})$. For convenience, we shall give the new settings.

For any positive integer $n$, an additive subgroup $G$ of $\Bbb{F}^n$ is called {\it nondegenerate} if $G$ contains an $\Bbb{F}$-basis of $\Bbb{F}^n$. Let $\ell_1,\;\ell_2$ and $\ell_3$ be three nonnegative integers such that
$$\ell=\ell_1+\ell_2+\ell_3>0.\eqno(3.3)$$
Take any nondegenerate additive subgroup $\G$ of $\Bbb{F}^{\ell_2+\ell_3}$ and $\G=\{0\}$ when $\ell_2+\ell_3=0$. Denote by $\Bbb{F}[t_1,t_2,...,t_{\ell_1+\ell_2}]$ the algebra of polynomials in $\ell_1+\ell_2$ variables over $\Bbb{F}$. Let ${\cal A}(\ell_1,\ell_2,\ell_3;\G)$ be a free $\Bbb{F}[t_1,t_2,...,t_{\ell_1+\ell_2}]$-module with the basis
$$\{x^{\al}\mid \al\in\G\}.\eqno(3.4)$$
Viewing ${\cal A}(\ell_1,\ell_2,\ell_3;\G)$ as a vector space over $\Bbb{F}$, we define a commutative associative algebraic operation ``$\cdot$'' on ${\cal A}(\ell_1,\ell_2,\ell_3;\G)$ by
$$(\zeta x^{\al})\cdot (\eta x^{\be})=\zeta\eta x^{\al+\be}\qquad\for\;\;\zeta,\eta\in \Bbb{F}[t_1,t_2,...,t_{\ell_1+\ell_2}],\;\al,\be\in\G.\eqno(3.5)$$
Note that  $x^{0}$ is the identity element, which is denoted as $1$ for convenience. When the context is clear, we shall omit the notion ``$\cdot$'' in any associative algebra product.

We define the linear transformations 
$$\{\ptl_{t_1},...,\ptl_{t_{\ell_1+\ell_2}},\ptl^{\ast}_1,...,\ptl^{\ast}_{\ell_2+\ell_3}\}\eqno(3.6)$$
on ${\cal A}(\ell_1,\ell_2,\ell_3;\G)$ by 
$$\ptl_{t_i}(\zeta x^{\al})=\ptl_{t_i}(\zeta)x^{\al},\;\;\ptl^{\ast}_j(\zeta x^{\al})=\al_j \zeta x^{\al}\eqno(3.7)$$
for $\zeta\in \Bbb{F}[t_1,t_2,...,t_{\ell_1+\ell_2}]$ and $\al=(\al_1,...,\al_{\ell_2+\ell_3})\in\G,$
where on $\Bbb{F}[t_1,t_2,...,t_{\ell_1+\ell_2}]$,  $\ptl_{t_i}$ are operators of taking partial derivative with respect to $t_i$. Then $\{\ptl_{t_1},...,\ptl_{t_{\ell_1+\ell_2}},\ptl^{\ast}_1,...,\ptl^{\ast}_{\ell_2+\ell_3}\}$ are mutually commutative derivations of ${\cal A}(\ell_1,\ell_2,\ell_3;\G)$. We set
$$\ptl_i=\ptl_{t_i},\;\;\ptl_{\ell_1+j}=\ptl^{\ast}_j+\ptl_{t_{\ell_1+j}},\;\;\ptl_{\ell_1+\ell_2+l}=\ptl_{\ell_2+l}^{\ast}\eqno(3.8)$$
for $i\in\ol{1,\ell_1},\;j\in\ol{1,\ell_2}$ and $l\in\ol{1,\ell_3}$. Then $\{\ptl_i\mid i\in\ol{1,\ell}\}$ is an $\Bbb{F}$-linearly independent set of derivations. 
Set
$${\cal D}=\sum_{i=1}^{\ell}\Bbb{F}\ptl_i\eqno(3.9)$$
and 
$${\cal W}(\ell_1,\ell_2,\ell_3;\G)= {\cal A}(\ell_1,\ell_2,\ell_3;\G){\cal D}.\eqno(3.10)$$
Then ${\cal W}(\ell_1,\ell_2,\ell_3;\G)$ is a standard form of the generalized simple Lie algebras of Witt type constructed in [X1]. Moreover, the Lie algebras found in [DZ1] are of the form ${\cal W}(\ell_1,0,\ell_3;\G)$.
In fact, we can rewrite the Lie algebra ${\cal W}$ in (3.1) under the condition (3.2) as ${\cal W}(\ell_1,\ell_2,\ell_3;\G)$ by considering the maximal $\Bbb{F}$-linearly independent subset of the set
 $$\{\{\al_i\mid (\al_1,...,\al_k)\in \G\}\mid i\in\ol{1,k}\}\eqno(3.11)$$
of $k$ sequences and changing variables in $\Bbb{F}[t_1,...,t_{k_1}]$ (cf. (2.1), (3.2)). 

Denote by $M_{m\times n}$ the set of $m\times n$ matrices with entries in $\Bbb{F}$ and by $GL_m$ the group of invertible $m\times m$ matrices. Set
$$G_{\ell_2,\ell_3}=\left\{\left(\begin{array}{cc}A&0_{\ell_2\times\ell_3}\\B& C\end{array}\right)\mid A\in GL_{\ell_2},\;B\in M_{\ell_2\times\ell_3},\;C\in  GL_{\ell_3}\right\},\eqno(3.12)$$
where $0_{\ell_2\times\ell_3}$ is the $\ell_2\times\ell_3$ matrix whose entries are zero. Then $G_{\ell_2,\ell_3}$ forms a subgroup of $GL_{\ell_2+\ell_3}$. Define an action of $G_{\ell_2,\ell_3}$ on $\Bbb{F}^{\ell_2+\ell_3}$ by
$$ g(\al)=\al g^{-1}\;\;(\mbox{matrix multiplication})\qquad\for\;\;\al\in\Bbb{F}^{\ell_2+\ell_3},\;g\in G_{\ell_2,\ell_3}.\eqno(3.13)$$
For any nondegenerate additive subgroup $\ups$ of $\Bbb{F}^{\ell_2+\ell_3}$ and $g\in G_{\ell_2,\ell_2}$, the set
$$g(\ups)=\{g(\al)\mid \al\in \ups\}\eqno(3.14)$$
also forms a nondegenerate additive subgroup of $\Bbb{F}^{\ell_2+\ell_3}$. Let
$$\Omega_{\ell_2+\ell_3}=\mbox{the set of nondegenerate additive subgroups of}\;\Bbb{F}^{\ell_2+\ell_3}.\eqno(3.15)$$
We have an action of $G_{\ell_2,\ell_3}$ on $\Omega_{\ell_2+\ell_3}$ by (3.14). Define the moduli space
$${\cal M}_{\ell_2,\ell_3}=\Omega_{\ell_2+\ell_3}/G_{\ell_2,\ell_3},\eqno(3.16)$$
which is the set of $G_{\ell_2,\ell_3}$-orbits in $\Omega_{\ell_2+\ell_3}$. 
\psp

{\bf Theorem 3.1}. {\it The Lie algebras} ${\cal W}(\ell_1,\ell_2,\ell_3;\G)$ {\it and} ${\cal W}(\ell_1',\ell_2',\ell_3';\G')$ {\it are isomorphic if and only if} $(\ell_1,\ell_2,\ell_3)=(\ell_1',\ell_2',\ell_3')$ {\it and there exists an element} $g\in  G_{\ell_2,\ell_3}$ {\it such that}
$g(\G)=\G'$. {\it In particular, there exists a one-to-one correspondence between the set of isomorphic classes of the Lie algebras of the form (3.10) and the following set:}
$$SW=\{(\ell_1,\ell_2,\ell_3,\varpi)\mid (0,0,0)\neq(\ell_1,\ell_2,\ell_3)\in\Bbb{N}^3,\;\varpi\in {\cal M}_{\ell_2,\ell_3}\}.\eqno(3.17)$$
{\it In other words, the set} $SW$ {\it is the structure space of the generalized simple Lie algebras of Witt type in the form (3.10)}. 
\psp

{\it Proof}. For convenience of proof, we redenote 
$$({\cal A}',{\cal D}',\ptl_i',t_j')\equiv ({\cal A},{\cal D},\ptl_i,t_j)\;\;\mbox{involved in}\;\;{\cal W}(\ell_1',\ell_2',\ell_3';\G').\eqno(3.18)$$
First, we assume $(\ell_1,\ell_2,\ell_3)=(\ell_1',\ell_2',\ell_3')$ and there exists an element $g\in G_{\ell_2,\ell_3}$  such that $g(\G)=\G'$. We write
$$g=\left(\begin{array}{cc}A& 0_{\ell_2\times\ell_3}\\B& C\end{array}\right)\qquad\mbox{with}\;\; A\in GL_{\ell_2},\;B\in M_{\ell_2\times\ell_3},\;C\in GL_{\ell_3}.\eqno(3.19)$$
By (3.13), we have
$$\G=\{\al'g\mid \al'\in\G'\}.\eqno(3.20)$$
Moreover, we set
$$\left(\begin{array}{c}\td{\ptl}'_{\ell_1+1}\\ \vdots\\ \td{\ptl}'_{\ell}\end{array}\right)=g^{-1}\left(\begin{array}{c}\ptl_{\ell_1+1}'\\ \vdots \\ \ptl'_{\ell}\end{array}\right),\;\;(\td{t}'_{\ell_1+1},...,\td{t}'_{\ell_1+\ell_2})=(t'_{\ell_1+1},...,t'_{\ell_1+\ell_2})A\eqno(3.21)$$
as matrix multiplications. For convenience, we let
$$\td{\ptl}_i'=\ptl_i',\;\;\td{t}_i'=t_i'\qquad\mbox{for}\;\;i\in\ol{1,\ell_1}.\eqno(3.22)$$
Now we define a linear map $\sgm$ from ${\cal W}(\ell_1,\ell_2,\ell_3;\G)$ to ${\cal W}(\ell_1,\ell_2,\ell_3;\G')$ by
$$\sgm(\prod_{i=1}^{\ell_1+\ell_2}t_i^{m_i}x^{\al}\ptl_j)=\prod_{i=1}^{\ell_1+\ell_2}(\td{t}'_i)^{m_i}x^{g(\al)}\td{\ptl}_j'\eqno(3.23)$$
for $(m_1,...,m_{\ell_1+\ell_2})\in\Bbb{N}^{\ell_1+\ell_2},\;\al\in\G$ and $j\in\ol{1,\ell}$. It is straightforward to verify that $\sgm$ is a Lie algebra isomorphism.

Next we assume that 
$$\sgm:{\cal W}(\ell_1,\ell_2,\ell_3;\G)\rightarrow {\cal W}(\ell_1',\ell_2',\ell_3';\G')\eqno(3.24)$$
is a Lie algebra isomorphism. We simply denote 
$${\cal W}={\cal W}(\ell_1,\ell_2,\ell_3;\G),\;\;{\cal W}'={\cal W}(\ell_1',\ell_2',\ell_3';\G')\eqno(3.25)$$
for convenience. Note that the adjoint operators 
$$\{\ad_d\mid d\in{\cal D}\},\;\;\{\ad_{d'}\mid d'\in{\cal D}'\}\;\;\mbox{are locally-finite}.\eqno(3.26)$$

Suppose that $\ell=\ell'=1$. If $\G\neq \{0\}$, we can prove that any element in ${\cal W}$ whose adjoint operator is locally-finite is in $\Bbb{F}\ptl_1$ by picking a well-order on $\Bbb{F}$ as an additive group (note $\Bbb{F}=\Bbb{C}$ is the most interesting case). The same statement holds for $(\G',W')$. A transformation $T$ of a vector space $V$ is called {\it locally-nilpotent} if for any $v\in V$, there exists a positive integer $m$ such that $T^m(v)=0$. Observe that $\ad_{\ptl_1}$ is locally-nilpotent if and only if  $\ad_{\sgm(\ptl_1)}$ is. Thus $\G=\{0\}$ if and only if $\G'=\{0\}$. If $\G=\G'=\{0\}$, we are done. Assume that $\G\neq\{0\}$ and $\G'\neq\{0\}$. By locally-finiteness,
$$\sgm{\ptl_1}=\lmd \ptl'_1\qquad\mbox{for some}\;\;0\neq\lmd\in\Bbb{F}.\eqno(3.27)$$
Note $\G,\G'\subset\Bbb{F}$ in this case. Moreover,
$$\Bbb{F}x^{\al}\ptl_1=\{u\in {\cal W}\mid [\ptl_1,u]=\al u\},\;\;\Bbb{F}x^{\al'}\ptl_1'=\{v\in {\cal W}\mid [\ptl_1',v]=\al' v\}\eqno(3.28)$$
for $\al\in\G$ and $\al'\in\G'$. Since
\begin{eqnarray*}\hspace{3cm}\al\sgm(x^{\al}\ptl_1)&=&\sgm([\ptl_1,x^{\al}\ptl_1])\\&=&
[\sgm(\ptl_1),\sgm(x^{\al}\ptl_1)]\\&=&\lmd [\ptl'_1,\sgm(x^{\al}\ptl_1)],\hspace{6.6cm}(3.29)\end{eqnarray*}
we have
$$\sgm(x^{\al}\ptl_1)=f(\al)x^{\lmd^{-1}\al}\ptl_1',\;\;\;f(\al)\in\Bbb{F}\qquad\for\;\;\al\in\G.\eqno(3.30)$$
Hence
$$\G'=\lmd^{-1} \G.\eqno(3.31)$$
Since $\ad_{\ptl_1}$ is semi-simple if and only $\ad_{\ptl_1'}$ is, we have
$$(\ell_1,\ell_2)=(\ell_1',\ell_2').\eqno(3.32)$$
Therefore, the theorem holds.

Next we assume $\ell,\ell'>1$. Set
$$X_{\al}=\Bbb{F}x^{\al},\;\;X'_{\al'}=\Bbb{F}x^{\al'},\;\;X=\sum_{\al\in\G}X_{\al},\;\;X'=\sum_{\al'\in\G'}X'_{\al'}.\eqno(3.33)$$
Then $X$ is the group algebra of $\G$ and $X'$ is the group algebra of $\G'$. Since ${\cal D}$ is locally-finite on ${\cal A}$, so is $\ad_{\cal D}$ on ${\cal W}$. This implies that $\ad_{\sgm({\cal D})}$ is locally-finite on ${\cal W}'$. Moreover, $\sgm({\cal D})$ is a commutative subalgebra of ${\cal W}'$. For any $d_1,d_2\in {\cal D}$ and $v\in{\cal A}'$, we have
$$[\sgm(d_1),v\sgm(d_2)]=\sgm(d_1)(v)\sgm(d_2).\eqno(3.34)$$
Hence $\sgm({\cal D})$ is locally-finite on ${\cal A}'$. For $i\in\ol{1,\ell_1}$, $\ptl_i$ is locally-nilpotent on ${\cal A}$. Naturally, $\ad_{\ptl_i}$ is locally-nilpotent on ${\cal W}$, which implies  
the  locally-nilpotency of $\ad_{\sgm(\ptl_i)}$ on ${\cal W}'$. Thus $\sgm(\ptl_i)$ is locally-nilpotent on ${\cal A}'$. By these facts, we have
$${\cal A}'=\bigoplus_{\be\in\Bbb{F}^{\ell_2+\ell_3}}\bar{\cal A}_{\be}\eqno(3.35)$$
with
$$\bar{\cal A}_{\be}=\{v\in{\cal A}'\mid (\sgm(\ptl_{\ell_1+j})-\be_j)^m(v)=0\;\mbox{for}\;j\in\ol{1,\ell_2+\ell_3}\;\mbox{and some}\;m\in\Bbb{N}\},\eqno(3.36)$$
where we have written $\be=(\be_1,...,\be_{\ell_2+\ell_3})$. Furthermore, we let
$$\bar{X}_{\be}=\{v\in{\cal A}'\mid \sgm(\ptl_i)(v)=0,\;\sgm(\ptl_{\ell_1+j})(v)=\be_j v\;\mbox{for}\;i\in\ol{1,\ell_1},\;j\in\ol{1,\ell_2+\ell_3}\}\eqno(3.37)$$
for $\be=(\be_1,...,\be_{\ell_2+\ell_3})\in\Bbb{F}^{\ell_2+\ell_3}.$ Obviously,
$$\bar{\cal A}_{\be}\neq\{0\}\Longleftrightarrow \bar{X}_{\be}\neq \{0\}.\eqno(3.38)$$

Set
$$\bar{\G}=\{\be\in\Bbb{F}^{\ell_2+\ell_3}\mid \bar{X}_{\be}\neq\{0\}\}.\eqno(3.39)$$
 Let $\al=(\al_1,...,\al_{\ell_2+\ell_3})\in \bar{\G}$ and let $0\neq z\in\bar{X}_{\al}$. For any $\ptl\in{\cal D}$, we have 
$$z\sgm(\ptl)=\sgm(w)\qquad\mbox{for some}\;\;w\in{\cal W}.\eqno(3.40)$$
Note 
$$0=[\sgm(\ptl_i),z\sgm(\ptl)]=\sgm([\ptl_i,w])\qquad\for\;\;i\in\ol{1,\ell_1}\eqno(3.41)$$
and
$$\al_j\sgm(w)=\al_jz\sgm(\ptl)=[\sgm(\ptl_{\ell_1+j}),z\sgm(\ptl)]=\sgm([\ptl_{\ell_1+j},w])\eqno(3.42)$$
for $j\in\ol{1,\ell_2+\ell_3}$. Hence
$$[\ptl_i,w]=0,\;\;[\ptl_{\ell_1+j},w]=\al_j w\qquad\for\;\;i\in\ol{1,\ell_1},\;j\in\ol{1,\ell_2+\ell_3}.\eqno(3.43)$$
Observe that
$$x^{\be}{\cal D}=\{u\in{\cal W}\mid [\ptl_i,u]=0,\;[\ptl_{\ell_1+j},u]=\be_j u\;\for\;i\in\ol{1,\ell_1},\;j\in\ol{1,\ell_2+\ell_3}\}\eqno(3.44)$$
for any $\be=(\be_1,...,\be_{\ell_2+\ell_3})\in \G$ and any root of ${\cal  W}$ with respect to $\ad_{\cal D}$ is in $\G$. Thus $\al\in \G$ and
$$w\in x^{\al}{\cal D}.\eqno(3.45)$$
So we obtain
$$\bar{\G}\subset \G.\eqno(3.46)$$
Moreover, we can write
$$z\sgm(\ptl)=\sgm(x^{\al}\tau_z(\ptl))\qquad\mbox{with}\;\;\tau_z(\ptl)\in{\cal D}.\eqno(3.47)$$
Hence we get an injective linear transformation $\tau_z$ on ${\cal D}$ because ${\cal A}'$ does not have zero divisors. Since ${\cal D}$ is finite-dimensional,
$\tau_z$ is a linear automorphism.  When $\al=0$, we get $0\neq z\in\Bbb{F}$ by (3.47). Therefore
$$\bar{X}_0=\Bbb{F}.\eqno(3.48)$$

Suppose $\G\neq\{0\}$. Then the elements of $\ad_{\cal D}$ are not all locally-nilpotent. So are those of  $\ad_{\sgm({\cal D})}$. Assume that $\bar{\G}=\{0\}$. Then $\sgm({\cal D})$ is locally-nilpotent. Since $\ad_{\sgm({\cal D})}$ is locally-finite, we can choose $\ptl\in {\cal D}$ and $d\in{\cal W}'$ such that
$$[\sgm(\ptl),d]=d.\eqno(3.49)$$
Since $d\neq 0$, there exists $v\in{\cal A}'$ such that 
$$d(v)\neq 0.\eqno(3.50)$$
On the other hand, there exists a positive integer $m$ such that
$$\sgm(\ptl)^m(v)=0.\eqno(3.51)$$
Note that
$$(\sgm(\ptl)-1)^{m+1}(d(v))=\sum_{j=0}^{m+1}(^{m+1}_{\;\;\;j})(\ad_{\sgm(\ptl)}-1)^j(d)\sgm(\ptl)^{m+1-j}(v)=0,
\eqno(3.52)$$
which contradicts the local nilpotency of $\sgm(\ptl)$. Thus $\bar{\G}\neq \{0\}$.

For any $\theta=(\theta_1,...,\theta_{\ell_2+\ell_3})\in\G$ and $\ptl=\sum_{i=1}^{\ell}a_i\ptl_i\in{\cal D}$ (cf. (3.3)), we define
$$\theta(\ptl)=\sum_{i=1}^{\ell_2+\ell_3}a_{\ell_1+i}\theta_i.\eqno(3.53)$$

Pick any $0\neq \al\in\bar{\G}$ and $0\neq z\in \bar{X}_{\al}$. We have (3.47). Let $0\neq \be\in\G$. Set
$$\gm=\be-\al.\eqno(3.54)$$
Since $\ell>1$, we can choose $0\neq \ptl\in{\cal D}$ such that
$$\gm(\ptl)=0.\eqno(3.55)$$
Moreover, we can pick  $d\in {\cal D}\setminus\Bbb{F}\tau_z(\ptl)$ such that 
$$\al(d)\neq 0.\eqno(3.56)$$
Then we have
\begin{eqnarray*}\hspace{3cm}[\sgm(x^{\gm}d), z\sgm(\ptl)]&=&\sgm(x^{\gm}d)(z)\sgm(\ptl)+z[\sgm(x^{\gm}d), \sgm(\ptl)]\\&=&\sgm(x^{\gm}d)(z)\sgm(\ptl)+z\sgm([x^{\gm}d,\ptl])\\&=&\sgm(x^{\gm}d)(z)\sgm(\ptl)-\gm(\ptl)z\sgm(x^{\gm}d)\\&=&\sgm(x^{\gm}d)(z)\sgm(\ptl),\hspace{5.5cm}(3.57)\end{eqnarray*}
by (3.55) and
\begin{eqnarray*}\hspace{3cm}[\sgm(x^{\gm}d),\sgm(x^{\al}\tau_z(\ptl))]&=&\sgm([x^{\gm}d,x^{\al}\tau_z(\ptl)])\\&=&\sgm(x^{\be}(\al(d)\tau_z(\ptl)-\gm(\tau_z(\ptl))d))\\ &\neq&0\hspace{7.1cm}(3.58)\end{eqnarray*}
by (3.56). Thus (3.47), (3.57) and (3.58) imply 
$$0\neq\sgm(x^{\gm}d)(z)\in \bar{X}_{\be}.\eqno(3.59)$$
Hence $\be\in\bar{\G}$. By (3.46), we obtain
$$\bar{\G}=\G.\eqno(3.60)$$

Now we want to prove that any nonzero element in $\bar{X}_{\al}$ is invertible for $\al\in\G$. Let $0\neq z\in \bar{X}_{\al}$. Pick any $0\neq z'\in \bar{X}_{-\al}$ and $\ptl\in {\cal D}$ such that $\al(\ptl)\neq 0$.
Since $\tau_z$ and $\tau_{z'}$ are invertible, we have
$$z\sgm(\tau_z^{-1}(\ptl))=\sgm(x^{\al}\ptl),\;\;z'\sgm(\tau_{z'}^{-1}(\ptl))=\sgm(x^{-\al}\ptl).\eqno(3.61)$$
Moreover,
$$[z'\sgm(\tau_{z'}^{-1}(\ptl)),z\sgm(\tau_z^{-1}(\ptl))]= zz'(\al(\tau_{z'}^{-1}(\ptl))\sgm(\tau_z^{-1}(\ptl))+\al(\tau_z^{-1}(\ptl))\sgm(\tau_{z'}^{-1}(\ptl))),\eqno(3.62)$$
$$[\sgm(x^{-\al}\ptl),\sgm(x^{\al}\ptl)]=\sgm([x^{-\al}\ptl,x^{\al}\ptl ])=2\al(\ptl)\sgm(\ptl).\eqno(3.63)$$
Hence
$$ zz'(\al(\tau_{z'}^{-1}(\ptl))\sgm(\tau_z^{-1}(\ptl))+\al(\tau_z^{-1}(\ptl))\sgm(\tau_{z'}^{-1}(\ptl)))=
2\al(\ptl)\sgm(\ptl)\neq 0\eqno(3.64)$$
by (3.61)-(3.63). Thus
$$0\neq zz'\in\Bbb{F}.\eqno(3.65)$$
This also shows 
$$\dim\bar{X}_{\al}=1\qquad\mbox{for}\;\;\al\in\G\eqno(3.66)$$
by (3.48). Furthermore, (3.35)-(3.37) imply
$${\cal A}'=\bar{X}\bar{\cal A}_{0}.\eqno(3.67)$$

Observe that 
$$(\bigcup_{\al\in\G}\bar{X}_{\al})\setminus \{0\}=\mbox{the set of all invertible elements in}\;{\cal A}'\eqno(3.68)$$
by (3.65)-(3.67). On the other hand, 
$$\mbox{the set of all invertible elements in}\;{\cal A}'=(\bigcup_{\al\in\G'}X'_{\al})\setminus \{0\}\eqno(3.69)$$
(cf. (3.33)) because ${\cal A}'={\cal A}(\ell_1',\ell_2',\ell_3';\G')$. Hence there exists a bijective map $\iota:\G\rightarrow \G'$ such that
$$\bar{X}_{\al}=X'_{\iota(\al)}\qquad\for\;\;\al\in\G.\eqno(3.70)$$
In particular,
$$\bar{X}=X'.\eqno(3.71)$$

Note that (3.71) implies
$$\bar{\cal A}_0=\Bbb{F}[t_1',t_2',...,t_{\ell'_1+\ell'_2}'].\eqno(3.72)$$
Set
$${\cal A}_{0,1}=\sum_{i=1}^{\ell_1+\ell_2}\Bbb{F}t_i,\;\;\bar{\cal A}_0^{(1)}=\{v\in \bar{\cal A}_0\mid \sgm({\cal D})(v)\subset\Bbb{F}\}.\eqno(3.73)$$

By the proof of Theorem 2.1 and the construction of ${\cal A}'={\cal A}(\ell_1',\ell_2',\ell_3';\G')$, we have
$$\dim (\bar{\cal A}_0^{(1)}/\Bbb{F})=\mbox{the transcendental degree of}\;\bar{\cal A}_0\;\mbox{over}\;\Bbb{F}=\ell'_1+\ell_2'.\eqno(3.74)$$
For any $z\in \bar{\cal A}^{(1)}_0$, we have
$$z\sgm(\ptl_1)=\sgm(w)\qquad\mbox{for some}\;\;w\in{\cal W}.\eqno(3.75)$$
Note
$$ \sgm(\ptl)(z)\sgm(\ptl_1)=[\sgm(\ptl),z\sgm(\ptl_1)]=\sgm([\ptl,w])\qquad\for\;\;\ptl\in{\cal D}.\eqno(3.76)$$
Since $\sgm(\ptl)(z)\in\Bbb{F}$, we obtain
$$[\ptl,w]=\sgm(\ptl)(z)\ptl_1\qquad\for\;\;\ptl\in{\cal D}.\eqno(3.77)$$
Since
$${\cal A}_{0,1}{\cal D}+{\cal D}=\{u\in {\cal W}\mid [{\cal D},u]\subset{\cal D}\}\eqno(3.78)$$ 
by the construction of ${\cal W}={\cal W}(\ell_1,\ell_3,\ell_3;\G)$, there exists a unique $\nu(z)\in  {\cal A}_{0,1}$ and $\rho(z)\in{\cal D}$ such that
$$w=\nu(z)\ptl_1+\rho(z)\;\;\mbox{and}\;\;\nu(z)=0 \Longleftrightarrow \sgm({\cal D})(z)=\{0\}\Longleftrightarrow z\in\Bbb{F}.\eqno(3.79)$$
Thus the map
$$z+\Bbb{F}\mapsto \nu(z) \eqno(3.80)$$
define an injective linear map from $\bar{\cal A}^{(1)}_0/\Bbb{F}$ to ${\cal A}_{0,1}$. So
$$\ell_1'+\ell_2'=\dim  (\bar{\cal A}_0^{(1)}/\Bbb{F})\leq \dim {\cal A}_{0,1}=\ell_1+\ell_2\eqno(3.81)$$
by (3.74). Exchanging positions of ${\cal W}$ and ${\cal W}'$, we can prove $ \ell_1+\ell_2\leq \ell_1'+\ell_2'$. Therefore,
$$ \ell_1+\ell_2=\ell_1'+\ell_2'.\eqno(3.82)$$

Denote
$$\hat{\ptl}_i=\sgm(\ptl_{\ell_1+i})|_{\bar{X}},\;\;\td{\ptl}_j=\ptl_{\ell_1'+j}'|_{\bar{X}}\qquad\for\;\;i\in\ol{1,\ell_2+\ell_3},\;j\in\ol{1,\ell'_2+\ell'_3}.\eqno(3.83)$$
Write
$$\sgm(\ptl_{\ell_1+i})=\sum_{j=1}^{\ell_1'}a_{i,j}\ptl'_j+\sum_{l=1}^{\ell'_2+\ell'_3}b_{i,l}\ptl'_{\ell_1'+l}\qquad\for\;\;i\in\ol{1,\ell_2+\ell_3},\eqno(3.84)$$
where
$$a_{i,j},b_{i,l}\in {\cal A}'\qquad \for\;\;i\in\ol{1,\ell_2+\ell_3},\;j\in\ol{1,\ell_1'},\;l\in\ol{1,\ell_2'+\ell_3'}.\eqno(3.85)$$
Since
$$\ptl'_j|_{\bar{X}}=0\qquad\for\;\;j\in\ol{1,\ell_1'},\eqno(3.86)$$
we have
$$\hat{\ptl}_i=\sum_{l=1}^{\ell'_2+\ell'_3}b_{i,l}\td{\ptl}_l\qquad\for\;\;i\in\ol{1,\ell_2+\ell_3}.\eqno(3.87)$$
Since $\G'$ is a nondegenerate additive subgroup of $\Bbb{F}^{\ell_2'+\ell_3'}$, there exists a basis $\Bbb{F}^{\ell_2'+\ell_3'}$:
$$\{\be^1,...,\be^{\ell_2'+\ell_3'}\}\subset \G'.\eqno(3.88)$$
Writing
$$\be^i=(\be^i_1,...,\be^i_{\ell'_2+\ell'_3})\qquad\for\;\;i\in\ol{1,\ell_2'+\ell_3'},\eqno(3.89)$$
we get an invertible matrix
$$\Psi=\left(\begin{array}{ccc}\be_1^1,&...,&\be_1^{\ell_2'+\ell_3'}\\ \vdots&\vdots&\vdots\\\be_{\ell_2'+\ell_3'}^1,&...,&\be_{\ell'_2+\ell_3'}^{\ell'_2+\ell_3'}\end{array}\right).\eqno(3.90)$$

Write
$$\iota^{-1}(\be^i)=(\al_1^i,...,\al_{\ell_2+\ell_3}^i)\in \G\qquad\for\;\;i\in\ol{1,\ell'_2+\ell_3'}.\eqno(3.91)$$
Note
$$\al_i^jx^{\be^j}=\hat{\ptl}_i(x^{\be^j})=\sum_{l=1}^{\ell'_2+\ell'_3}b_{i,l}\td{\ptl}_l(x^{\be^j})=(\sum_{l=1}^{\ell'_2+\ell'_3}b_{i,l}\be^j_l)x^{\be^j},\eqno(3.92)$$
equivalently,
$$\al_i^j=\sum_{l=1}^{\ell'_2+\ell'_3}b_{i,l}\be^j_l\eqno(3.93)$$
for $i\in\ol{1,\ell_2+\ell_3}$ and $j\in\ol{1,\ell_2'+\ell_3'}$. Thus
$$(\al_i^1,...,\al_i^{\ell_2'+\ell_3'})= (b_{i,1},...,b_{i,\ell_2'+\ell_3'})\Psi\qquad\for\;\;i\in\ol{1,\ell_2+\ell_3}.\eqno(3.94)$$
In particular,
$$(b_{i,1},...,b_{i,\ell_2'+\ell_3'})=(\al_i^1,...,\al_i^{\ell_2'+\ell_3'})\Psi^{-1}\in\Bbb{F}^{\ell_2'+\ell_3'}\qquad\for\;\;i\in\ol{1,\ell_2+\ell_3}.\eqno(3.95)$$
Since $\{\hat{\ptl}_1,...,\hat{\ptl}_{\ell_2+\ell_3}\}$ is linearly independent due to the nondegeneracy of $\bar{\G}=\G$ and $\{\td{\ptl}_1,...,\td{\ptl}_{\ell_2'+\ell_3'}\}$ is linearly independent because of the nondegeneracy of $\G'$,
we have
$$\ell_2+\ell_3\leq \ell_2'+\ell_3'\eqno(3.96)$$
by (3.87) and (3.95). Exchanging positions of ${\cal W}$ and ${\cal W}'$, we can prove $\ell_2'+\ell_3'\leq \ell_2+\ell_3$. Hence we have 
$$\ell_2+\ell_3=\ell_2'+\ell_3'.\eqno(3.97)$$
Thus the matrix
$$B=(b_{i,j})_{(\ell_2+\ell_3)\times(\ell_2+\ell_3)}\eqno(3.98)$$
is nondegenerate. 

Observe that 
$$\sgm(\ptl_{\ell_1+\ell_2+i})\;\mbox{is semi-simple because}\;\ad_{\sgm(\ptl_{\ell_1+\ell_2+i})}\;\mbox{is}\eqno(3.99)$$
for $i\in\ol{1,\ell_3}$. Thus
$$0=\sgm(\ptl_{\ell_1+\ell_2+i})(t'_{\ell_1'+j})=b_{\ell_2+i,j}\qquad\for\;\;i\in\ol{1,\ell_3},\;j\in\ol{1,\ell_2'}\eqno(3.100)$$
by (3.72) and (3.84). So
$$\hat{\ptl}_{\ell_2+i}=\sum_{j=1}^{\ell_3'}b_{\ell_2+i,\ell_2'+j}\td{\ptl}_{\ell_2'+j}\qquad\for\;\;i\in\ol{1,\ell_3}.\eqno(3.101)$$
The above expression implies $\ell_3\leq \ell_3'$. Exchanging positions of ${\cal W}$ and ${\cal W}'$, we can similarly prove $\ell_3'\leq \ell_3$. Thus
$$\ell_3=\ell_3'.\eqno(3.102)$$
Therefore, we obtain
$$(\ell_1,\ell_2,\ell_3)=(\ell_1',\ell_2',\ell_3')\eqno(3.103)$$
by (3.82), (3.97) and (3.102).

Let $\al=(\al_1,...,\al_{\ell_2+\ell_3})\in \G$. Denote $\be=(\be_1,...,\be_{\ell_2+\ell_3})=\iota(\al)$. We have
$$\al_ix^{\be}=\sgm(\ptl_{\ell_1+i})(x^{\be})=(\sum_{j=1}^{\ell_2+\ell_3}b_{i,j}\ptl'_{\ell_1+j})(x^{\be})=(\sum_{j=1}^{\ell_2+\ell_3}b_{i,j}\be_j)(x^{\be}),\eqno(3.104)$$
equivalently
$$\al_i=\sum_{j=1}^{\ell_2+\ell_3}b_{i,j}\be_j\eqno(3.105)$$
for $i\in\ol{1,\ell_2+\ell_3}$. Hence we get
$$\iota(\al)=\be=\al (B^t)^{-1},\eqno(3.106)$$
where $B^t$ is the transpose of $B$. By (3.95) and (3.100), we have $B^t\in  G_{\ell_2,\ell_3}$ (cf. (3.12)). Denoting $g=B^t$, we get
$$\G'=g(\G)\eqno(3.107)$$
by (3.13) and (3.106). This completes the proof of Theorem 3.1.

\vspace{1cm}

\noindent{\Large \bf References}

\hspace{0.5cm}

\begin{description}

\item[{[DZ1]}] D. \v{Z}. Dokovick and K. Zhao, Generalized Cartan type W Lie algebras in characteristic 0, {\it J. Algebra} {\bf 195} (1997), 170-210.

\item[{[DZ2]}] D. \v{Z}. Dokovick and K. Zhao, Derivarions, isomorphisms, and second cohomology of generalized Witt algebras, {\it Trans. Amer. Math. Soc.} {\bf 350} (1998), 643-664. 

\item[{[Ka1]}] V. G. Kac, A description of filtered Lie algebras whose associated graded Lie algebras are of Cartan types, {\it Math. of USSR-Izvestijia} {\bf 8} (1974), 801-835.

\item[{[Ka2]}] V. G. Kac, Lie superalgebras, {\it Adv. Math.} {\bf 26} (1977), 8-96.

\item[{[K]}] N. Kawamoto, Generalizations of Witt algebras over a field of characteristic zero,
{\it Hiroshima Math. J.} {\bf 16} (1986), 417-462.

\item[{[O]}] J. Marshall Osborn, New simple infinite-dimensional Lie algebras of characteristic 0, {\it J. Algebra} {\bf 185} (1996), 820-835.

\item[{[P]}] D. P. Passman, Simple Lie algebras of Witt type, {\it J. Algebra} {\bf 206} (1998), 682-692.

\item[{[X1]}] X. Xu, New generalized simple Lie algebras of Cartan type over a field with characteristic 0, {\it J. Algebra}, in press.

\item[{[X2]}] X. Xu, {\it Algebraic Theory of Hamiltonian Superoperators}, monograph, to appear.

\item[{[Z]}] K. Zhao, Isomorphisms between generalized Cartan type W Lie algebras in characteristic zero, {\it Canadian J. Math.} {\bf 50} (1998), 210-224.

\end{description}
\end{document}